\documentclass[12pt,a4paper,oneside]{amsart}
\usepackage{fullpage}
\usepackage{amssymb,amsfonts}
\usepackage{epsfig}
\usepackage{graphicx}

\usepackage{graphics}
\newtheorem{theorem}{Theorem}[section]
\newtheorem{lemma}[theorem]{Lemma}

\theoremstyle{definition}
\newtheorem{definition}[theorem]{Definition}
\newtheorem{remark}[theorem]{Remark}

\numberwithin{equation}{section}

\DeclareMathOperator{\diam}{diam} 

\newcommand{\be}{\begin{equation}}
\newcommand{\ee}{\end{equation}}

\DeclareMathOperator{\rad}{rad}

\makeindex

\def\Xint#1{\mathchoice 
 {\XXint\displaystyle\textstyle{#1}}%
{\XXint\textstyle\scriptstyle{#1}}%
{\XXint\scriptstyle\scriptscriptstyle{#1}}%
 {\XXint\scriptscriptstyle\scriptscriptstyle{#1}}%
 \!\int}
\def\XXint#1#2#3{{\setbox0=\hbox{$#1{#2#3}{\int}$}
 \vcenter{\hbox{$#2#3$}}\kern-.5\wd0}}

 \def\dashint{\Xint-}

\begin{document}
\title{Generalized Lebesgue points for Sobolev functions}
\author{Nijjwal Karak}
\address{Department of Mathematics and Statistics, University of Jyv\"askyl\"a, P.O. Box 35, FI-40014, Jyv\"askyl\"a, Finland}
\email{nijjwal.n.karak@jyu.fi}
\thanks{This work was supported by the Academy of Finland via the Centre of Excellence in Analysis and Dynamics Research (Grant no. 271983)}
\begin{abstract}
In this article, we show that a function $f\in M^{s,p}(X),$ $0<s\leq 1,$ $0<p<1,$ where $X$ is a doubling metric measure space, has generalized Lebesgue points outside a set of $\mathcal{H}^h$-Hausdorff measure zero for a suitable gauge function $h.$
\end{abstract}
\maketitle
\section{Introduction}
By the Lebesgue differentiation theorem, almost every point in $\mathbb{R}^n$ is a Lebesgue point of a locally integrable function, that is
\begin{equation*}
\lim_{r\rightarrow 0}\frac{1}{\vert B(x,r)\vert}\int_{B(x,r)}u(y)\,dy=u(x)
\end{equation*}
for almost every $x\in\mathbb{R}^n$ and for a locally integrable function $u.$ It is a well-known fact that a function $f\in W^{1,p}(\mathbb{R}^n),$ $1\leq p\leq n,$ has Lebesgue points outside a set of $p$-capacity zero, \cite{EG92}, \cite{Zie89}, \cite{HKM06}. Recently, there has been some interests in studying Lebesgue points for Sobolev functions on metric measure spaces, specially for functions in Haj\l asz-Sobolev space $M^{1,p}(X)$ and in Newtonian space (or Sobolev space) $N^{1,p}(X)$ defined by Haj\l asz \cite{Haj96} and Shanmugalingam \cite{Sha00} respectively. The usual argument for obtaining the existence of Lebesgue points outside a small set for a Sobolev function goes as follows. First of all, Lebesgue points exist outside a set of capacity zero, see \cite{KL02}, \cite{KKST08} for Sobolev functions on metric measure spaces. Secondly, each set of positive Hausdorff $h$-measure, for a suitable $h,$ is of positive capacity, see \cite{KM72}, \cite{AH96}, \cite{Oro87} for sets in $\mathbb{R}^n$ and \cite{BO05}, \cite{KK15A}, \cite{KKST08} for sets in metric measure spaces. Combining these results one gets the existence of Lebesgue points outside a set of Hausdorff $h$-measure zero, for a suitable $h,$ see \cite{KK15B} for more details on this.\\

In this paper, we study the existence of Lebesgue points of a function in Haj\l asz-Sobolev space $M^{s,p}(X),$ for $0<s\leq 1,$ $0<p<1,$ outside a small set in terms of Hausdorff $h$-measure. Recall that A measurable function $f : X\rightarrow\mathbb{R},$ where $X=(X,d,\mu)$ is a metric measure space, belongs to the Haj\l asz-Sobolev space $M^{s,p}(X),$ $0<s\leq 1,$ $p>0,$ if and only if $f\in L^p(X)$ and there exists a nonnegative function $g\in L^p(X)$ such that the inequality
\begin{equation}\label{Hajlasz}
\vert f(x)-f(y)\vert\leq d(x,y)^s(g(x)+g(y))
\end{equation}
holds for all $x,y\in X\setminus E,$ where $\mu(E)=0.$ This definition is due to Haj\l asz for $s=1,$ \cite{Haj96} and to Yang for fractional scales, \cite{Yan03}.\\

Recently, Heikkinen, Koskela and Tuominen have studied the existence of generalized Lebesgue points for functions in $M^{s,p}(X),$ $0<s\leq 1,$ $0<p<\infty,$ outside a set of capacity zero \cite{HKT}. They have also studied the same for functions in Haj\l asz-Besov spaces $N^s_{p,q}$ and Haj\l asz-Triebel-Lizorkin spaces $M^s_{p,q}.$ Notice that $M^s_{p,\infty}(X)=M^{s,p}(X),$ see \cite{KYZ11}. The existence of Lebesgue points outside a small set in terms of capacity for Besov and Triebel-Lizorkin functions has been studied in \cite{AH96}, \cite{HN07}, \cite{Net89} and the relation between Besov-capacity and Hausdorff measure has been studied in \cite{Cos09}. In this paper we only consider functions in Haj\l asz-Sobolev spaces and we avoid the use of capacity here. Here we use medians to define generalized Lebesgue points, as we do not have the integrability of the functions. Median allows us to study the oscillation of measurable functions. Please see Section 2 for the definitions of medians and generalized Lebesgue points. Medians have been studied for example in \cite{PT12}, \cite{Fuj91}, \cite{FZ72}, \cite{Str79}.\\
Our result is the following.
\begin{theorem}\label{main}
Let $(X,d,\mu)$ be a doubling metric measure space. Let $f\in M^{s,p}(X),$ where $0<s\leq 1,$ $0<p<1.$ Then $\lim_{r\rightarrow 0}\,m_f^{\gamma}(B(z,r))$ exists
outside a set $E_{\epsilon}$ with $\mathcal{H}^h(E_\epsilon)=0,$ whenever
$$h(B(x,\rho))=\frac{\mu(B(x,\rho))}{\rho^{sp}}\log^{-p-\epsilon}(1/\rho)$$
for any $\epsilon>0.$
\end{theorem}
We refer to Section 2 for the definition of generalized Hausdorff $h$-measure and also for the existence of the above limit outside a small set for a measurable and finite almost everywhere function.\\
   
Note that this result is new even in $\mathbb{R}^n.$ For $f\in M^{1,p}(\mathbb{R}^n),$ where $\frac{n}{n+1}<p<1,$ we can use integral averages instead of medians by the recent result of Koskela and Saksman \cite{KS08}. They have proved that if a function $f\in M^{1,p}(\mathbb{R}^n),$ for $p>\frac{n}{n+1},$ then $f$ is locally integrable.\\

The paper is organized as follows. We explain our notations and state a couple of elementary results in Section 2. The proof of the theorem is given in Section 3 and we also give a couple of remarks there.\\

\indent {\textit{Acknowledgement}.} I wish to thank Professor Pekka Koskela for his fruitful suggestions.
\section{Notation and preliminaries}
We assume throughout that $X=(X,d,\mu)$ is a metric measure space equipped with a metric $d$ and a Borel regular outer measure $\mu.$ We call such a $\mu$ as a measure. The Borel-regularity of the measure $\mu$ means that all Borel sets are $\mu$-measurable and that for every set $A\subset X$ there is a Borel set $D$ such that $A\subset D$ and $\mu(A)=\mu(D).$\\

We denote open balls in $X$ with center $x\in X$ and radius $0<r<\infty$ by $$B(x,r)=\{y\in X : d(y,x)<r\}.$$
If $B=B(x,r)$ is a ball and $\lambda>0,$ we write
$$\lambda B=B(x,\lambda r).$$
With small abuse of notation we write $\rad(B)$ for the radius of a ball $B$ and we always have $$\diam(B)\leq 2\rad(B),$$
and the inequality can be strict.\\

A Borel regular measure $\mu$ on a metric space $(X,d)$ is called a \textit{doubling measure} if every ball in $X$ has positive and finite measure and there exist a constant $C_{\mu}\geq 1$ such that
\begin{equation*}
\mu(B(x,2r))\leq C_{\mu}\,\mu(B(x,r))
\end{equation*}
for each $x\in X$ and $r>0.$ We call a triple $(X,d,\mu)$ a \textit{doubling metric measure space} if $\mu$ is a doubling measure on $X.$\\

If $A\subset X$ is a $\mu$-measurable set with finite and positive measure, then the \textit{integral average} of a function $u\in L^1(A)$ over $A$ is
$$u_A=\dashint_A u\, d\mu=\frac{1}{\mu(A)}\int_A u\, d\mu.$$
\\
\begin{definition}
Let $0<\gamma\leq 1/2.$ The $\gamma$-median $m_f^{\gamma}(A)$ of a measurable, almost everywhere finite function $f$ over a set $A\subset X$ of finite measure is
 \begin{equation*}
m_f^{\gamma}(A)=\max\{M\in\mathbb{R}:\mu(\{x\in A:f(x)<M\})\leq\gamma\mu(A)\}.
 \end{equation*}
\end{definition}
We mention here two basic properties of medians, for the proof see \cite{PT12}, \cite{HKT}.\\
(i) If $f$ is continuous, then for every $x\in X$ and $0<\gamma\leq 1/2,$
\begin{equation*}
\lim_{r\rightarrow 0}\,m_f^{\gamma}(B(x,r))=f(x).
\end{equation*}
(ii) There exists a set $E$ with $\mu(E)=0$ such that
\begin{equation*}
\lim_{r\rightarrow 0}\,m_f^{\gamma}(B(x,r))=f(x)
\end{equation*}
holds for every $0<\gamma\leq 1/2$ and $x\in X\setminus E.$\\
\begin{definition}
Let $0<\gamma\leq 1/2$ and $f$ be a measurable, almost everywhere finite function. A point $x\in X$ is a  \textit{generalized Lebesgue point} of $f,$ if
\begin{equation*}
\lim_{r\rightarrow 0}\,m_f^{\gamma}(B(x,r))=f(x).
\end{equation*}
\end{definition}
We recall that the \textit{generalized Hausdorff $h$-measure} is defined by
\begin{equation*}
\mathcal{H}^h(E)=\limsup_{\delta\rightarrow 0}H_{\delta}^{h}(E),
\end{equation*}
where
\begin{equation*}
H_{\delta}^{h}(E)=\inf\left\{\sum h(B(x_i,r_i)) : E\subset\bigcup B(x_i,r_i),~r_i\leq\delta\right\},
\end{equation*}
where the dimension gauge function $h$ is required to be continuous and increasing with $h(0)=0,$ see \cite{KKST08}.\\
\indent Given a non-negative, locally integrable function $f$ on $\mathbb{R}^n,$ its Riesz potential is defined as
\begin{equation}\label{riesz1}
I_{\alpha}f(x)=C(\alpha, n)\int_{\mathbb{R}^n}\frac{f(y)}{\vert x-y\vert^{n-\alpha}}\,dy,
\end{equation}
or its local version
\begin{equation}\label{riesz2}
I_{\alpha}^{\Omega}f(x)=C(\alpha, n)\int_{\Omega}\frac{f(y)}{\vert x-y\vert^{n-\alpha}}\,dy,
\end{equation}
where $0<\alpha<n$ and $C(\alpha, n)$ is a suitable constant.\\
\indent For the convenience of reader we state here a fundamental covering lemma (for a proof see \cite[2.8.4-6]{Fed69} or \cite[Theorem 1.3.1]{Zie89}).
\begin{lemma}[5B-covering lemma]\label{cover}
Every family $\mathcal{F}$ of balls of uniformly bounded diameter in a metric space $X$ contains a pairwise disjoint subfamily $\mathcal{G}$ such that for every $B\in\mathcal{F}$ there exists $B'\in\mathcal{G}$ with $B\cap B'\neq\emptyset$ and $\diam(B)<2\diam(B').$ In particular, we have that
$$\bigcup_{B\in\mathcal{F}}B\subset\bigcup_{B\in\mathcal{G}}5B.$$
\end{lemma}

\section{Proof of Theorem \ref{main}}

\begin{proof}
Fix $\epsilon>0$. Let $h$ be as in the statement of the Theorem. Let us write $B_j=B(z,2^{-j})$ for $z\in X$ and $j\in\mathbb{N}.$ Our first aim is to show that the sequence $(m_f^{\gamma}(B_j))_j$ is a Cauchy sequence outside a set of $\mathcal{H}^h$-measure zero. From the definition of median it easily follows that, for all $j\in\mathbb{N},$
\begin{equation}\label{median_1}
\mu(B_{j}^l)=\mu\big(\big\{x\in B_{j}:f(x)\leq m_f^{\gamma}(B_{j})\big\}\big)\geq(1-\gamma)\mu(B_{j})
\end{equation}
and
\begin{equation}\label{median_2}
\mu(B_{j}^u)=\mu\big(\big\{y\in B_j:f(y)\geq m_f^{\gamma}(B_j)\big\}\big)\geq\gamma\mu(B_j).
\end{equation}
Then by using the inequality \eqref{Hajlasz} and the Fubini theorem, we obtain
\begin{align*}
\begin{split}
\mu(B_j^u)\,\mu(B_{j+1}^l)\,\vert m_f^{\gamma}(B_j)-m_f^{\gamma}(B_{j+1})\vert^p &\leq \int_{B_j^u}\int_{B_{j+1}^l}\vert f(x)-f(y)\vert ^p\,d\mu(x)\,d\mu(y)\\
&\leq \int_{B_j^u}\int_{B_{j+1}^l}d(x,y)^{sp}\big(g(x)+g(y)\big)^p\,d\mu(x)\,d\mu(y)\\
&\leq 2^p\int_{B_j^u}\int_{B_{j+1}^l}d(x,y)^{sp}\big(g^p(x)+g^p(y)\big)\,d\mu(x)\,d\mu(y)\\
&\leq 2^{2p}2^{-spj}\int_{B_j^u}\int_{B_{j+1}^l}g^p(x)\,d\mu(x)\,d\mu(y)\\
&\qquad+ 2^{2p}2^{-spj}\int_{B_j^u}\int_{B_{j+1}^l}g^p(y)\,d\mu(x)\,d\mu(y)\\
&= 2^{2p}2^{-spj}\mu(B_{j}^u)\int_{B_{j+1}^l}g^p(x)\,d\mu(x)\\
&\qquad+ 2^{2p}2^{-spj}\mu(B_{j+1}^l)\int_{B_j^u}g^p(x)\,d\mu(x).
\end{split}
\end{align*}
Using the doubling property and the inequalities \eqref{median_1} and \eqref{median_2}, we get 
\begin{eqnarray*}
\vert m_f^{\gamma}(B_j)-m_f^{\gamma}(B_{j+1})\vert ^p &\leq & \frac{2^{2p}2^{-spj}}{\mu(B_{j+1}^l)}\int_{B_{j+1}^l}g^p(x)\,d\mu(x)+ \frac{2^{2p}2^{-spj}}{\mu(B_{j}^u)}\int_{B_j^u}g^p(x)\,d\mu(x)\\
&= & 2^{2p}2^{-spj}\bigg[\frac{\mu(B_j)}{\mu(B_{j+1}^l)}+\frac{\mu(B_j)}{\mu(B_{j}^u)}\bigg]\dashint_{B_j}g^p(x)\,d\mu(x)\\
&\leq & 2^{2p}2^{-spj}\bigg[\frac{C_{\mu}}{1-\gamma}+\frac{1}{\gamma}\bigg]\dashint_{B_j}g^p(x)\,d\mu(x)\\
&= & C2^{-spj}\dashint_{B_j}g^p(x)\,d\mu(x),
\end{eqnarray*}
where $C=C(\gamma,p,C_{\mu}).$
For $m,l\in\mathbb{R}^n,$ $m<l,$ let us consider the difference 
\begin{eqnarray}\label{difference}
\nonumber \vert m_f^{\gamma}(B_l)-m_f^{\gamma}(B_m)\vert &\leq &\sum_{j=m}^{l-1}\vert m_f^{\gamma}(B_j)-m_f^{\gamma}(B_{j+1})\vert\\
&\leq & C\sum_{j=m}^{l-1}2^{-sj}\left(\dashint_{B_j}g^p(x)\,d\mu(x)\right)^\frac{1}{p}.
\end{eqnarray}
\indent Let $h_1(B(x,\rho))=\frac{\mu(B(x,\rho))}{\rho^{sp}}\log^{-p-\epsilon/2}(1/\rho).$ If we have $\int_{B(z,r)}g^p\,dx\leq Ch_1(B(z,r))$ for all sufficiently small $0<r<1/5,$ then $(m_f^{\gamma}(B_j))_j$ is a Cauchy sequence, by \eqref{difference}. On the other hand, let us consider the set
\begin{multline*}
E_\epsilon=\bigg\{z\in\mathbb{R}^n :\text{there exists arbitrarily small}~0<r_z<\frac{1}{5}~\text{such that}\\ \int_{B(z,r_z)}g^p\, d\mu(x)\geq Ch_1(B(z,r_z))\bigg\}.
\end{multline*}
Let $0<\delta<1/5.$ Then we get a pairwise disjoint family $\mathcal{G}$ consisting of balls as above, by using the 5B-covering lemma, such that
$$E_{\epsilon}\subset\bigcup_{B\in\mathcal{G}}5B,$$
where $\diam(B)<2\delta$ for $B\in\mathcal{G}.$ Then we estimate
\begin{eqnarray*}
\mathcal{H}_{10\delta}^{h_1}(E_{\epsilon}) & \leq & C\sum_{B\in\mathcal{G}} h_1\left(B(z,\rad(B))\right)\\
& \leq & C\sum_{B\in\mathcal{G}} \int_{B}g^p\, d\mu(x)\\
& \leq & C\int_{\bigcup\limits_{B\in\mathcal{G}}B}g^p\, d\mu(x)<\infty.
\end{eqnarray*}
It follows that $\mathcal{H}^{h_1}(E_{\epsilon})<\infty$ and hence we have that $\mathcal{H}^{h}(E_{\epsilon})=0,$ which gives us the existence of $ \lim_{j\rightarrow\infty}\,m_f^{\gamma}(B(z,2^{-j}))$ for $\mathcal{H}^{h}$-a.e. $z\in X.$\\
\indent For given $r>0,$ we can always find $j\in\mathbb{N}$ such that $2^{-(j+1)}<r<2^{-j}.$ By using the same method as above we conclude that
\begin{equation*}
\vert m_f^{\gamma}(B_j)-m_f^{\gamma}(B(z,r))\vert\leq C2^{-spj}\dashint_{B_j}g^p(x)\,d\mu(x)
\end{equation*}
and that $\lim_{r\rightarrow 0}\,m_f^{\gamma}(B(z,r))$ exists outside $E_{\epsilon}.$
\end{proof}
\begin{remark}
It is known that $f\in M^{1,1}(X)$ has Lebesgue points outside a set $E$ with $\mathcal{H}^{h}(E)=0$ with $h(B(x,\rho))=\frac{\mu(B(x,\rho))}{\rho}$ provided $X$ supports a $1$-Poincar\'e inequality, \cite{KKST08}. We do not know if one can obtain a better result than Theorem \ref{main} for $f\in M^{1,p}(X)$ by showing that the exceptional set has $\mathcal{H}^{h}$-Hausdorff measure zero with $h(B(x,\rho))=\frac{\mu(B(x,\rho))}{\rho^p}.$ In $\mathbb{R}^n,$ one possible approach is to use a Riesz potential after the inequality \eqref{difference}, as shown below.\\
It is easy to see, from \eqref{difference}, that
\begin{eqnarray*}
\vert m_f^{\gamma}(B_l)-m_f^{\gamma}(B_m)\vert &\leq & C\left(\sum_{j=m}^{l-1}2^{-jp}\dashint_{B_j}g^p(x)\,dx\right)^\frac{1}{p}\\
&\leq & C\left(\int_{B_m}\frac{g^p(x)}{\vert z-x\vert^{n-p}}\,dx\right)^\frac{1}{p}\\
&=& CI_{p}^{B_m}g^p(z).
\end{eqnarray*}
Then we use Theorem 3.1.4 (a) of \cite{AH96} to conclude that $\lim_{r\rightarrow 0}\,m_f^{\gamma}(B(z,r))$ exists outside $E$ with $\mathcal{L}^n(E)=0.$ It would be interesting to know if there is a similar estimate as in Theorem 3.1.4 (a) of \cite{AH96} for the $\mathcal{H}^{n-\alpha}$-Hausdorff measure of the set $\{z:I_{\alpha}u(z)>\lambda\},$ for $u\in L^1(\mathbb{R}^n),$ $0<\alpha<n$ and for all $\lambda>0.$ This would improve our result in this case. 
\end{remark}
\begin{remark}
In $\mathbb{R}^n,$ for the case when $n/(n+1)<p<1,$ we use telescoping arguments between the centred balls and also use inequality \eqref{Hajlasz} to get similar estimate as in \eqref{difference} for the integral averages instead of medians. Similar technique can be found in \cite{HK00}. Then it is easy to see that $\lim_{r\rightarrow 0}\,f_{B(z,r)}$ exists outside a set of $\mathcal{H}^h$-measure zero with the same $h$ as in Theorem \ref{main}.
\end{remark}
\def\bibname{References}
\bibliography{lebesgue_p<1}
\bibliographystyle{alpha}
\end{document}